\documentclass[11pt]{article}
\title{Fock Space Decomposition of L\'{e}vy Processes}
\author{R. F. Streater,\\Department of Mathematics,\\King's College
London,\\Strand, London,\\WC2R 2LS}
\date{10/1/2002}

\newtheorem{definition}{Definition}

\begin{document}
\maketitle \setlength{\oddsidemargin}{0in}
\setlength{\evensidemargin}{0in}
\begin{abstract}
We show that the general Levy process can be embedded in a suitable 
Fock space, classified by cocycles of the real line regarded as a group.
The formula of de Finetti corresponds to coboundaries. Kolmogorov's 
processes correspond to cocycles of which the derivative is
a cocycle of the Lie algebra of ${\bf R}$. The L\'{e}vy formula
gives the general cocycle.
\end{abstract}

\section{Cyclic representations of groups}
Let $G$ be a group and $g\mapsto U_g$ a multiplier cyclic
representation of $G$ on a Hilbert space ${\cal H}$, with
multiplier $\sigma:G\times G\rightarrow {\bf C}$ and cyclic vector
$\Psi$. This means that
\begin{itemize}
\item $U_gU_h=\sigma(g,h)U_{gh}$ for all $g,h\in G$.
\item $u(e)=I$ where $e$ is the identity of then group and $I$ is
the identity operator on ${\cal H}$.
\item Span$\,\left\{U_g\Psi:g\in G\right\}$ is dense in ${\cal
H}$.
\end{itemize}
If $\sigma=1$ we say that $U$ is a true representation.

Recall that a multiplier of a group $G$ is a measurable
two-cocycle in $Z^2(G,U(1))$; so $\sigma$ is a map $G\times
G\rightarrow U(1)$ such that $\sigma(e,g)=\sigma(g,e)=1$ and
\begin{equation}
\sigma(g,h)\sigma(g,hk)^{-1}\sigma(gh,k)\sigma(h,k)^{-1}=1.
\label{sigma}
\end{equation}
(\ref{sigma}) expresses the associativity of operator
multiplication of the $U(g)$. $\sigma$ is a coboundary if there is
a map $b:G\rightarrow U(1)$ with $b(e)=1$ and
\begin{equation}
\sigma(g,h)=b(gh)/(b(g)b(h)).
\label{2boundary}
\end{equation}
\vspace{.1in} \noindent We also need the concept of a one cocycle
$\psi$ in a Hilbert space ${\cal K}$ carrying a unitary
representation $V$. $\psi$ is a map $G\rightarrow {\cal K}$
such that
\begin{equation}
V(g)\psi(h)=\psi(gh)-\psi(g) \hspace{.5in}\mbox{for }g,h\in G.
\end{equation}
$\psi$ is a coboundary if there is a vector $\psi_0\in{\cal K}$
such that
\begin{equation}
\psi(g)=(V(g)-I)\psi_0.
\label{oneboundary}
\end{equation}
Coboundaries are always cocycles. We say that, in
(\ref{2boundary}) and (\ref{oneboundary}), $\sigma$ is the
coboundary of $b$ and $\psi$ is the coboundary of $\psi_0$.

\newpage
We say that two cyclic  $\sigma$-representations $\{{\cal
H},U,\Psi\}$ and $\{{\cal K},V,\Phi\}$ are {\em cyclically
equivalent} if there exists a unitary operator $W:{\cal
H}\rightarrow{\cal K}$ such that $V_g=WU_gW^{-1}$ for all $g\in
G$, and
$W\Psi=\Phi$. Any cyclic multiplier representation $\{{\cal
H},U,\Psi\}$ defines a function $F$ on the group by
\begin{equation}
F(g):=\langle\Psi,U_g\Psi\rangle,
\end{equation}
which satisfies $\sigma$-positivity:
\begin{eqnarray}
F(e)&=&1\\
\sum_{ij}\overline{\lambda}_i\lambda_j\sigma(g_i^{-1},g_j)F(g_i^{-1}g_j)
&\geq& 0.
\label{twisted}
\end{eqnarray}
$F$ is called the characteristic function of the representation,
because
\begin{itemize}
\item Two cyclic multiplier representations of $G$ are cyclically
equivalent if and only if they have the same characteristic
function;
\item Given a function on $G$ satisfying $\sigma$-positivity, then there
exists a cyclic $\sigma$-representation of which it is the
characteristic function.
\end{itemize}
\noindent If $G=\{s\in{\bf R}\}$, $\sigma=1$ and $U_s$ is
continuous, then $F$ obeys the hypotheses of Bochner's theorem and
defines a probability measure $\mu$ on ${\bf R}$. More generally,
we can apply Bochner's theorem (if $\sigma=1$) to any
one-parameter subgroup $s\mapsto g(s)\in G_0\subseteq\in G$. Then
$U_{g(s)},\;s\in{\bf R}$ is a one-parameter unitary group; its
infinitesimal generator is a self-adjoint operator $X$ on ${\cal
H}$. The relation to $\mu$ is given as follows: let $X=\int\lambda
dE(\lambda)$ be the spectral resolution of $X$. Then
\begin{equation}
\mu(\lambda_1,\lambda_2]=\langle\Psi,\left(E(\lambda_2)-E(\lambda_1)\right)
\Psi\rangle.
\end{equation}
Conversely, given any random variable $X$ on a probability space
$(\Omega,\mu)$, we can define the cyclic unitary representation of
the group ${\bf R}$ by the multiplication operator
\begin{equation}
U(s)=\exp\{isX\}
\end{equation}
and use the cyclic vector $\Psi(\omega)=1$ on the Hilbert space
$L^2(\Omega, d\mu)$. In this way, probability theory is reduced to
the study of cyclic representations of abelian groups, and quantum
probability to the study of cyclic $\sigma$-representations of
non-abelian groups.

\section{Processes as Tensor Products}
Given a cyclic $\sigma$-representation $\{{\cal H},U,\Phi\}$ of a
group $G$, we can get a multiplier representation of the product
group $G^n:=G\times G\times \ldots \times G$ ($n$ factors) on
${\cal H}\otimes{\cal H}\ldots\otimes{\cal H}$, by acting on the
vector $\Psi\otimes\Psi\ldots\otimes\Psi$ by the unitary operators
$U(g_1,\ldots,g_n):=U(g_1)\otimes\ldots U(g_n)$, as each $g_j$
runs over the group $G$. The resulting cyclic $\sigma^{\otimes
n}$-representation is denoted
\begin{equation}
\left\{{\cal H}^{\otimes n},U^{\otimes n},\Psi^{\otimes n}\right\}.
\end{equation}
The twisted positive function on $G^n$ defined by this cyclic
representation is easily computed to be
\begin{equation}
F^{\otimes n}(g_1,\ldots,g_n)=F(g_1)F(g_2)\ldots F(g_n).
\end{equation}
If $G$ has a one-parameter subgroup $G_0$, then the infinitesimal
generators $X_j$ of this subgroup in the $j^{\rm th}$ place define
random variables $(j=1,\ldots, n)$ that are all independent in the measure
$\mu^{\otimes n}$ on ${\bf R}^n$ defined by $F^{\otimes n}$, and are
all identically distributed. They can thus be taken as the increments
of a process in discrete time $t=1,\ldots,n$. To get a process with time
going to infinity, we can embed each tensor product ${\cal H}^{\otimes n}$
in the ``incomplete infinite tensor product'' of von Neumann, denoted
\begin{equation}
\bigotimes_{j=1}^{j=\infty}\rule{.5cm}{0cm}^\Psi{\cal H}_j\hspace{.5in}\mbox{where }
{\cal H}_j={\cal H}\mbox{ for all }j.
\end{equation}

It is harder to construct processes in
continuous time. We made \cite{RFS} the following definition:
\begin{definition}
{\em A cyclic $G$-representation $\{{\cal H},U,\Psi\}$ is said to
be {\em infinitely divisible} if for each positive integer $n$
there exists another cyclic $G$-representation $\{{\cal
K},V,\Phi\}$ such that $\{{\cal H},U,\Psi\}$ is cyclically
equivalent to $\{{\cal K}^{\otimes n},V^{\otimes n}, \Phi^{\otimes
n}\}$.}
\end{definition}
The picturesque notation $\left\{{\cal H}^{\otimes \frac{1}{n}},
U^{\otimes\frac{1}{n}},\Psi^{\otimes\frac{1}{n}}\right\}$ can be
used for $\{{\cal K},V,\Phi\}$.

If $G={\bf R}$ then $\{{\cal H},U,\Psi\}$
is infinitely divisible if and only if the corresponding measure
$\mu$ given by Bochner's theorem is infinitely divisible
\cite{RFS}. It is clear that $\{{\cal H},U,\Psi\}$ is infinitely
divisible if and only if there exists a branch of
$F(g)^{\frac{1}{n}}$ which is positive semi-definite on $G$.

This criterion was extended in \cite{PS} to
$\sigma$-representations. In that case, for each $n$, there should
exist an $n^{\rm th}$ root $\sigma(g,h)^{\frac{1}{n}}$ which is
also a multiplier. One can then consider cyclic representations
such that for each $n$, $F(g)^{\frac{1}{n}}$ has a branch which is
$\sigma^{\frac{1}{n}}$-positive semi-definite.

If $\{{\cal H},U,\Psi\}$ is an infinitely
divisible $G$-representation, then we may construct a continuous
tensor product of the Hilbert spaces ${\cal H}_t$, where $t\in{\bf
R}$ and all the Hilbert spaces are the same. This gives us, in the
non-abelian case, quantum stochastic processes with independent
increments. See references. The possible constructions are
classified in terms of cocycles of the group $G$. Here we shall
limit discussion to the analysis of the L\'{e}vy formula in these
terms.
\section{The cocycle}

Let $F: G\rightarrow {\bf C}$ and $F(e)=1$. It is a classical result
for $G={\bf R}$ that a function $F^{\frac{1}{n}}$
has a branch that is positive semidefinite for all $n>0$ if and only if
$\log F$ has a branch $f$ such that $f(0)=0$ and $f$ is {\em
conditionally} positive semidefinite. This is equivalent to
$f(x-y)-f(x)-f(-y)$ being
positive semidefinite. This result is easily extended to groups
\cite{RFS} and $\sigma$-representations \cite{PS,PS2}. Let us
consider the case where $\sigma=1$. It follows that an infinitely divisible
true cyclic representation $\{{\cal H},U,\Psi\}$ of $G$ defines a
conditionally positive semidefinite function $f(g)={\rm log}\langle\Psi,
U(g)\Psi\rangle$, so that
\begin{equation}
\sum_{j,k}\overline{\alpha}_j\alpha_k\left(\rule{0cm}{1cm}f\left(
g_j^{-1}g_k\right)-f(g_j)^{-1}-f(g_k)\right)\geq 0.
\label{qform}
\end{equation}
We can use this positive semidefinite form to make Span$\,G$ into a
pre-scalar product space, by defining
\begin{equation}
\langle\psi(g),\psi(h)\rangle:=f(g^{-1}h)-f(g^{-1})-f(h),
\hspace{.5in}g,h\in G.
\label{scalar}
\end{equation}
Let ${\cal K}$ be the Hilbert space, that is the separated and
completed space got this way. There is a natural injection
$\psi:G\rightarrow{\cal K}$, namely, $g\mapsto [g]$, the
equivalence class of $g$ given by the relation $g\sim h$ if the
seminorm defined by (\ref{qform}) vanishes on $g-h$. The left
action of the group $G$ on this function is not quite unitary; in
fact the following is a unitary representation \cite{A}:
\begin{equation}
V(h)\psi(g):=\psi(hg)-\psi(h).
\end{equation}
One just has to check from (\ref{scalar}) that the group law
$V(g)V(h)=V(gh)$ holds, and that
\begin{equation}
\langle V(h)\psi(g_1),V(h)\psi(g_2)\rangle=\langle\psi(g_1),\psi(g_2)\rangle.
\end{equation}
Thus we see that $\psi(g)$ is a one-cocycle relative to the
$G$-representation $V$ \cite{A}.
\section{The embedding theorem}
Given a Hilbert space ${\cal K}$, the {\em Fock space} defined by ${\cal K}$
is the direct sum of all symmetric tensor products of ${\cal K}$,
\begin{equation}
EXP\,{\cal K}:={\bf C}\bigoplus{\cal K}\bigoplus({\cal
K}\otimes{\cal K})_s \bigoplus\ldots .
\end{equation}
The element $1\in{\bf C}$ is called the Fock vacuum.
The following {\em coherent states} form a total set in $EXP\,{\cal K}$:
\begin{equation}
EXP\,\psi:=1+\psi+(1/2!)\psi\otimes\psi+\ldots,
\hspace{.5in}\psi\in{\cal K}.
\end{equation}
The notation is natural, in view of the easy identity
\begin{equation}
\langle
EXP\,\psi(g),EXP\,\psi(h)h\rangle=\exp\{\langle\psi(g),\psi(h)\rangle\}.
\end{equation}
Then the embedding theorem \cite{RFS} says that if $\{{\cal
H},U,\Psi\}$ is an infinitely divisible cyclic representation,
then it is cyclically equivalent to the cyclic representation $W$
on $EXP\,{\cal K}$, with the Fock vacuum as the cyclic vector,
with the unitary representation $W(h)$defined on the total set of
coherent states by
\begin{equation}
W(h)EXP\,\psi(g)=F(hg)/F(g)EXP\,\psi(hg).
\end{equation}
The proof is simply a verification. This result has been called
\cite{RFS,PS} the Araki-Woods embedding theorem; more properly
this name belongs to the embedding \cite{RFS} of the {\em process}
that one constructs from $\{{\cal H},U,\Phi\}$, which is similar
to a deep result in \cite{AW}. For the group ${|bf R}$ it amounts
to the Wiener chaos expansion.

\section{The L\'{e}vy Formula}
Although the above theory was developed for quantum mechanics, it
includes the theory of L\'{e}vy processes, which is the class of
processes with independent increments. This is just the case when
the group in question is ${\bf R}$, or ${\bf R}^n$. The latter
group has projective representations in $n>1$, and using these
leads to the free quantised field \cite{RFS0}. The true
representations lead to generalised random fields.

Every projective representation of ${\bf R}$ is a true
representation, which is multiplicity free if it is cyclic. By
reduction theory, it is then determined by a measure on the dual
group, here ${\bf R}$. Araki showed that a one-cocycle can be
algebraic or topological. For the group ${\bf R}$, the algebraic
cocycles are all of the form $f(x-y)-f(x) -f(-y)=axy$. This is
satisfied by the Gaussian term $log F(x)=-\frac{a}{2}x^2+ibx$, and
this is the only possibility. The Poisson($\lambda$) is an example
of a coboundary, when $\log F(t) =c\lambda(e^{ipt}-1)$ for some
$p$, the increment of the jumps. The weighted mixture of these
coboundaries gives di Finetti's formula \cite{D}:
\begin{equation}
\log F(t)=\lambda\left\{ibt-\frac{a^2t^2}{2}+c\int\left(e^{ipt}-1\right)
dP(p)\right\}.
\end{equation}
That this is not the most general infinitely divisible measure was
recognised by Kolmogorov \cite{K}. In our terms, this is the statement that
not all topological cocycles are coboundaries (the cohomology is
non-trivial). Kolmogorov considered random variables with finite variance
relative to the measure $dP$. This is equivalent in our terms to
$dP=|\widehat{\psi}(p)|^2dp$ and the
cocyle $\psi$ being of the form $\psi(x)=(V(x)-I)\psi_0$, where
$i\partial_x\psi_0$ is square integrable over the group ${\bf R}$, but
$\psi_0$ might not be. Thus,
$\psi$ is a cocycle for the Lie algebra of the group, a case treated in
\cite{RFS3}. This gives us Kolmogorov's formula
\begin{eqnarray}
\log F(t)&=&\lambda\left\{ibt-a^2t^2/2+\right.\nonumber\\
&+&\left.\int\left(e^{ipt} -1-itp\right)|\psi(p)|^2dp\right\}.
\end{eqnarray}
The term $\int(-itp)|\psi(p)|^2dp$ is possibly divergent near
$p=0$ but is not required to exist on its own near $p=0$, since
the function $M=e^{ipt}-1-ipt$ behaves as $p^2$ near the origin.
But to retain a meaning, Kolmogorov's formula does need
$p|\widehat{\psi}(p)|^2$ to be integrable at infinity. This is not
needed for the general cocycle, so the formula is not the most
general.

\noindent L\'{e}vy gave the answer \cite{L} by replacing $M$ by
\begin{equation}
e^{ipt}-1-ipt/(1+p^2),
\end{equation}
so that $\widehat{\psi}$ has no constraint at infinity other than
being $L^2$. The general form of an infinitely divisible random
variable, the L\'{e}vy formula, in effect constructs the most
general cocycle of the group ${\bf R}$ by requiring only that
$p\widehat{\psi}(p)$ should be locally square-integrable at $p=0$,
and $\widehat{\psi}(p)$ should be square-integrable at all other
points.


\begin{thebibliography}{99}
\bibitem{A} Araki, H., Factorizable representations of current
algebra, non-commutative extension of the L\'{e}vy-Khinchin
formula, and cohomology of a solvable group with values in a
Hilbert space,  Publ. Res Inst. Math Sci (RIMS) Kyoto, 1970/71.
\bibitem{AW} Araki, H., and Woods, E. J., Complete Boolean Lattices
of Type I von Neumann Algebras, Publ. Res. Inst. Math. Sci., (Kyoto),
{\bf 2}, 157- ,1966.
\bibitem{D} de Finetti, B., Sulla funzione a incremento aleatorio, Atti
Acad. Naz. Lincei. Rend. Cl. Sci. Fis. Mat. Nat. {\bf 10}, 163-168,
325-329, 548-553, 1929.
\bibitem{EF} Erven, J., and Falkowski, B.-J., Low order cohomology
and applications, Lecture Notes in Maths., {\bf 877},
Springer-Verlag, 1981.
\bibitem{F} Falkowski, F.-J., Factorizable and infinitely divisible
PUA representations of locally compact groups, J. of Mathematical Phys.,
{\bf 15}, 1060-1066, 1974.
\bibitem{G} Guichardet, A., Symmetric Hilbert spaces and related
topics, Lecture Notes in Maths., {\bf 261}, Springer-Verlag, 1972.
\bibitem{K} Kolmogorov, A. N., Sulla forma generale di un processo
stocastico omogeneo, Atti Acad. Naz. Lincei Rend. Cl. Sci. Fis. Mat.
Nat. {\bf 15}, 805-808, 866-869 1932.
\bibitem{L} L\'{e}vy, P., Sur les int\'{e}grales dont les \'{e}l\'{e}ments
sont des variables al\'{e}atoires ind\'{e}pendentes, Annali R.
Scuola Norm. Sup. Pisa, {\bf 3}, 337-366, 1934; {\bf 4}, 217-218, 1935.
\bibitem{PS} Parthasarathy, P. K., and Schmidt, K., Positive definite
kernels, continuous tensor products, and central limit theorems...,
Lecture Notes in Maths., {\bf 272}, Springer-Verlag, 1972.
\bibitem{PS2} Parthasarathy, K. R., and Schmidt, K., Factorizable
representations of current groups and the Araki-Woods embedding
theorem, Acta Math., {\bf 128}, 53-71, 1972.
\bibitem{RFS0} Streater, R. F. Current Commutation Relations and
Continuous Tensor Products, Nuovo Cimento, {\bf 53A}, 487-495, 1968.
\bibitem{RFS} Streater, R. F., Current commutation relations, continuous
tensor products, and infinitely divisible group representations;
pp 247-263, in {\em Local Quantum Theory} ed R. Jost, Academic
Press, 1969.
\bibitem{RFS2} Streater, R. F., A continuum analogue of the lattice gas,
Commun. Math. Phys., {\bf 12}, 226-232, 1969.
\bibitem{RFS3} Streater, R. F., Infinitely divisible
representations of Lie algebras, Zeits. fur Wahr. und verw. Gebiete,
{\bf 19}, 67-80, 1971.
\bibitem{VGG} Vershik, A. M., Gelfand, I. M., and Graev, M. J.,
Irreducible representations of the group $G^X$ and cohomologies,
Funct. Anal. and Appl. (translated from Russian), {\bf 8}, 67-69,
1974.
\end{thebibliography}
\end{document}